\documentclass[psamsfonts, reqno, 12pt]{amsart}

\usepackage{amsmath,amsfonts,amsthm,amssymb,psfrag,graphicx}

\DeclareMathOperator\ord{ord}

\numberwithin{equation}{section}

\newtheorem{thm}{Theorem}[section]
\newtheorem{lem}[thm]{Lemma}
\newtheorem{cor}[thm]{Corollary}
\newtheorem{prop}[thm]{Proposition}

\theoremstyle{definition}
\newtheorem{eg}[thm]{Example}
\newtheorem{definition}[thm]{Definition}

\newtheorem{rmk}[thm]{Remark}

\newcommand{\fix}{F}
\newcommand{\rat}{\mathsf{Q}}
\newcommand{\sph}{\mathsf{S}}
\newcommand{\nes}{\mathsf{N}}
\newcommand{\nesn}{\mathsf{N}^{\mathsf{n}}}
\newcommand{\nesv}{\mathsf{N}^{\mathsf{v}}}
\newcommand{\ch}{\textrm{char}}
\newcommand{\ass}{\textrm{Asc}}

\newcommand{\coht}{\textrm{coht}}
\newcommand{\hght}{\textrm{ht}}
\newcommand{\kdim}{\textrm{kdim}}
\newcommand{\exprk}{\textrm{exprk}}
\newcommand{\entrk}{\textrm{entrk}}

\title[Periodic point data detects subdynamics]{Periodic Point
Data Detects Subdynamics in Entropy Rank One}
\author{Richard Miles}
\author{Thomas Ward}
\date{\today}

\address{School of Mathematics, University of East Anglia, Norwich, NR4 7TJ, UK}
\subjclass[2000]{22D40, 37A15, 37A35}
\thanks{This research was supported by E.P.S.R.C. grant EP/C015754/1.}

\begin{document}

\begin{abstract}
A framework for understanding the geometry of continuous actions
of~$\mathbb Z^d$ was developed by Boyle and Lind using the notion of
expansive behavior along lower-dimensional subspaces. For
algebraic~$\mathbb Z^d$-actions of entropy rank one, the expansive
subdynamics is readily described in terms of Lyapunov exponents.
Here we show that periodic point counts for elements of an entropy rank
one action determine the expansive subdynamics. Moreover, the finer
structure of the non-expansive set is visible in the topological and
smooth structure of a set of functions associated to the periodic
point data.
\end{abstract}

\maketitle

\section{Introduction}

Let~$\beta$ be an action of~$\mathbb Z^d$ by homeomorphisms of a
compact metric space~$(X,\rho)$; thus for each~$\mathbf n\in\mathbb
Z^d$ there is an associated homeomorphism~$\beta^{\mathbf n}$,
and~$\beta^{\mathbf m}\circ\beta^{\mathbf n}=\beta^{\mathbf
m+\mathbf n}$ for all~$\mathbf m,\mathbf n\in\mathbb Z^d$. Such an
action is called \emph{expansive} if there is some~$\delta>0$ with
the property that if~$x,y$ are distinct points in~$X$ then there is
some~$\mathbf n$ for which~$\rho(\beta^{\mathbf n}x,\beta^{\mathbf
n}y)>\delta$. Any such~$\delta$ is called an expansive constant for
the action. Boyle and Lind~\cite{MR1355295} introduced the following
notion, which reveals a rich geometrical structure inside an
expansive action. A subset~$A\subset\mathbb R^d$ is called
\emph{expansive for~$\beta$}, or~$\beta$ is \emph{expansive
along~$A$}, if there exist constants~$\delta>0$ and~$t>0$ with the
property that
\[
\sup_{\mathbf n,d(\mathbf n,A)<t}\rho(\beta^{\mathbf
n}x,\beta^{\mathbf n}y)\le\delta\implies x=y\mbox{ for all }x,y\in X
\]
where~$d(\mathbf n,A)$ denotes the distance from the point~$\mathbf
n$ to the set~$A$ in the Euclidean metric on~$\mathbb R^d$. Of
particular importance is the behavior along subspaces.
Write~$\mathsf G_k$ for the Grassmannian of~$k$-dimensional
subspaces of~$\mathbb R^d$; this is a compact~$k(d-k)$-dimensional
manifold in the usual topology (subspaces are close if their
intersections with the unit~$(d-1)$-sphere~$\sph_{d-1}$ are
close in the Hausdorff topology). Following Boyle and Lind, write
\[
\nes_k(\beta)=\{V\in\mathsf G_k\mid V\mbox{ is not expansive
for~$\beta$}\}.
\]
The main structural result from~\cite{MR1355295} is that if~$X$ is
infinite, then~$\nes_{d-1}(\beta)$ is a non-empty compact set, and
the set~$\nes_{d-1}(\beta)$ governs all of the non-expansive
behavior in the sense that any element of~$\nes_k(\beta)$ must be a
subspace of some element of~$\nes_{d-1}(\beta)$. For algebraic
systems, in which~$X$ is a compact metric group and each
map~$\beta^{\mathbf n}$ is a continuous group automorphism, the
subdynamical structure was determined by Einsiedler, Lind, Miles and
Ward~\cite{MR1869066}, where a finer structure was found inside the
set~$\nes_{d-1}(\beta)$ reflecting the two different ways in which
an algebraic dynamical system can fail to be expansive.

A different insight into a topological~$\mathbb Z^d$ action is a
combinatorial one coming from periodic points. Write~$\fix_{\mathbf
n}(\beta)=\{x\in X\mid\beta^{\mathbf n}x=x\}$ for the set of points
fixed by the homeomorphism~$\beta^{\mathbf n}$. The combinatorial
data of all these numbers may be thought of as a map
\[
\mathbf n\mapsto\vert\fix_{\mathbf n}(\beta)\vert\in\mathbb
N\cup\{\infty\},
\]
where~$\infty$ denotes the cardinality of an infinite compact group.
Our purpose here is to show that the combinatorial data contained in
this map determines the expansive subdynamics for a certain class of
systems. These systems are the expansive algebraic systems of
entropy rank one, which in particular means that the
set~$\fix_{\mathbf n}(\beta)$ is finite except in degenerate
situations.

\section{Ranks and Subdynamics}

The following notions come from~\cite[Sect.~7]{MR1869066}.
Let~$\beta$ be an action of~$\mathbb Z^d$ by homeomorphisms of a
compact metric space~$(X,\rho)$ as before. The \emph{expansive rank}
of~$\beta$ is
\[
\exprk(\beta)=\min\{k\mid\nes_{k}(\beta)\neq\mathsf G_k\},
\]
that is the smallest dimension in which some expansive subspaces are
seen. The \emph{entropy rank} of~$\beta$ is
\[
\entrk(\beta)=\max\{k\mid\mbox{there is a rational~$k$-plane~$V$
with }h(\beta,V)>0\},
\]
where~$h(\beta,V)$ denotes the topological entropy of the~$\mathbb
Z^{\dim(V)}$-action given by restricting~$\beta$ to~$V\cap\mathbb
Z^d$. By~\cite[Prop.~7.2]{MR1869066},
\[
\entrk(\beta)\le\exprk(\beta).
\]

Algebraic~$\mathbb Z^d$-actions have a convenient description in
terms of commutative algebra due to Kitchens and
Schmidt~\cite{MR1036904} which we will need.
Let~$R_d=\mathbb{Z}[u_1^{\pm 1}, \dots, u_d^{\pm 1}]$ be the ring of
Laurent polynomials in commuting variables~$u_1, \dots, u_d$ with
integer coefficients. If~$X$ is a compact metrizable abelian group
and~$\alpha$ is a~$\mathbb{Z}^d$-action by continuous
automorphisms~$\alpha^{\mathbf{n}}$ of~$X$, then the Pontryagin dual
group~$M=\widehat{X}$ has the structure of a discrete
countable~$R_d$-module, obtained by first identifying the dual
automorphism~$\widehat{\alpha}^{\mathbf{n}}$ with multiplication by
the monomial~$u^{\mathbf{n}}=u_1^{n_1}\dots u_d^{n_d}$, and then
extending additively to multiplication by polynomials. Conversely,
for any countable~$R_d$-module~$M$, there is an
associated~$\mathbb{Z}^d$-action on a compact group obtained by
dualizing the action induced by multiplying by monomials on~$M$. A
full account of this correspondence and the resulting theory is
given in Schmidt's monograph~\cite{MR1345152}. An important aspect
of this approach is the interpretation of dynamical properties as
algebraic properties of~$M$, particularly in terms of the set of
associated prime ideals of~$M$, written~$\ass(M)$. We will describe
systems as Noetherian if they correspond to Noetherian modules, and
in the reverse direction will describe modules as having various
dynamical properties if the corresponding system has those
properties.

The simplest algebraic systems are those corresponding to cyclic
modules~$R_d/\mathfrak p$ for a prime ideal~$\mathfrak p\subset
R_d$, and these will be called \emph{prime actions}. This gives a
third natural notion of `rank' to an algebraic~$\mathbb Z^d$-action.
Recall that the \emph{Krull dimension}~$\kdim(S)$ of a ring~$S$ is
the maximum of the lengths~$r$ taken over all strictly decreasing
chains~$\mathfrak p_0\supset\mathfrak
p_1\supset\cdots\supset\mathfrak p_r$ of prime ideals in~$S$ (see
Matsumura~\cite[Chap.~1\S5]{MR1011461}). Boyle and
Lind~\cite[Th.~7.5]{MR1355295} show that if~$\mathfrak p$ is a prime
ideal generated by~$g$ elements, then
\[
\exprk(\alpha_{R_d/\mathfrak p})\ge\kdim(R_d/\mathfrak p)\ge d-g
\]
and
\[
\exprk(\alpha_{R_d/\mathfrak p})\ge d-g+1.
\]
Moreover,~\cite[Prop.~7.3]{MR1869066} shows that
\[
\entrk(\alpha_{R_d/\mathfrak p})=\kdim(R_d/\mathfrak
p)\le\exprk(\alpha_{R_d/\mathfrak p})
\]
if~$\mathfrak p$ is non-principal. The
\emph{height}~$\hght(\mathfrak p)$ of a prime ideal~$\mathfrak p
\subset R_d$ is equal to the Krull dimension of~$R_d$ localized at
$\mathfrak p$, equivalently the maximal length~$r$ of a strictly
decreasing chain of prime ideals
\[
\mathfrak p=\mathfrak p_0\supset \mathfrak p_1\supset
\cdots\supset \mathfrak p_r=(0).
\]
The \emph{co-height}~$\coht(\mathfrak p)$ of~$\mathfrak p$ is equal
to the Krull dimension of the domain~$R_d/\mathfrak p$, equivalently
it is the maximal length~$r$ of a strictly increasing chain of prime
ideals
\[
\mathfrak p=\mathfrak p_0\subset\mathfrak
p_1\subset\cdots\subset\mathfrak p_r.
\]
The domain~$R_d$ is universally catenary~\cite[Prop.~18.9,
Cor.~18.10]{MR1322960} and hence~\cite[Th.~13.8]{MR1322960} shows
that for each~$\mathfrak p\in R_d$,
\[
\hght(\mathfrak p)+\coht(\mathfrak p)=\kdim(R_d)=d+1.
\]

Using associated primes, Einsiedler and Lind~\cite{MR2031042}
provide the following classification of entropy rank one actions for
which the associated module~$M$ is Noetherian (see
Proposition~\ref{main_structure_theorem}). When~$M$ is not
Noetherian, problems arise in relation to finding the set of
possible entropy values for general
algebraic~$\mathbb{Z}^d$-actions; this is closely related to
Lehmer's problem and is discussed more fully in~\cite{MR2031042}.

\begin{prop}\label{main_structure_theorem}
Let~$\alpha_M$ be a Noetherian algebraic~$\mathbb{Z}^d$-action. Then
\begin{enumerate}
\item $\alpha_M$ has entropy rank one if and only if each of the
associated prime actions~$\alpha_{R_d/\mathfrak{p}}$ has entropy
rank one. Equivalently, for each
prime~$\mathfrak{p}\in\ass(M)$,~$\textrm{coht}(\mathfrak{p})\leq1$;
\item $\alpha_M$ has expansive rank one if and only if each of the
associated prime actions~$\alpha_{R_d/\mathfrak{p}}$ has expansive
rank one; and
\item if~$\alpha_M$ is expansive then~$\alpha_M$ has expansive rank one if
and only if~$\alpha_M$ has entropy rank one.
\end{enumerate}
\end{prop}

\begin{proof}
See~\cite[Prop.~4.4 and~6.1, Th.~7.1 and~7.2]{MR2031042}.
\end{proof}

In particular, an expansive rank one action may also be thought of
as an expansive entropy rank one action. Further properties of
entropy rank one actions are discussed in~\cite{MR2031042}
and~\cite{miles_zeta}. Of particular importance is the observation
that if~$\textrm{coht}(\mathfrak{p})=1$ then the field of
fractions~$K$ of the domain~$R_d/\mathfrak{p}$ is a \emph{global
field} by~\cite[Prop.~6.1]{MR2031042}. Moreover, the places of~$K$,
denoted by~$\mathcal{P}(K)$, are determined by the
ideal~$\mathfrak{p}$. From this infinite set of places, we isolate
\[
\mathcal{S}_\mathfrak{p}=\{w\in\mathcal{P}(K)\mid w \textrm{ is
unbounded on } R_d/\mathfrak{p}\}.
\]
Here~$w$ being \emph{unbounded} means that~$|R_d/\mathfrak{p}|_w$ is
an unbounded subset of~$\mathbb{R}$. Note
that~$\mathcal{S}_\mathfrak{p}$ contains all the infinite places
of~$K$. Furthermore,~$\mathcal{S}_\mathfrak{p}$ is finite
because~$R_d/\mathfrak{p}$ is finitely generated.

The description of expansive subdynamics for algebraic~$\mathbb
Z^d$-actions is further refined in~\cite[Sect.~8]{MR1869066} to
reflect the two ways in which an algebraic dynamical system can fail
to be expansive. It can fail in a way which relates to the
Noetherian condition for modules, and this failure will result in a
set of directions denoted~$\nesn$. It can also fail to be expansive
in the way a quasihyperbolic toral automorphisms fails to be
expansive, by having the higher-rank analogue of an eigenvalue with
unit modulus; this failure arising from the varieties of the
associated prime ideals results in a set of directions
denoted~$\nesv$.

The Noetherian condition is described as follows.
Each~$\mathbf{n}\in\mathbb{Z}^d$ defines a half-space
\[
H=\{\mathbf x\in\mathbb R^d\mid\mathbf x\cdot\mathbf n\le
0\}\subset\mathbb R^d,
\]
which has an associated ring~$R_H=\mathbb Z[u^{\mathbf m}\mid\mathbf
m\in H\cap\mathbb Z^d]\subset R_d$. A module over~$R_d$ is also a
module over~$R_H$.

\begin{definition}
Let~$M$ be a Noetherian~$R_d$-module and let~$V\subset\mathbb R^d$
be a~$k$-dimensional subspace. Then~$M$ is said to be
\emph{Noetherian along~$V$} if~$M$ is a Noetherian~$R_H$-module for
every half-space~$H$ containining~$V$. The collection of
all~$k$-dimensional subspaces along which~$M$ is not Noetherian is
denoted~$\nesn_k(\alpha_M)$.
\end{definition}

For the variety condition, let~$\mathfrak a\subset R_d$ be any
ideal. Write
\[
\mathsf V(\mathfrak a)=\{\mathbf z\in(\mathbb C\setminus\{0\})^d\mid
f(\mathbf z)=0\mbox{ for all }f\in R_d\}
\]
and define the \emph{amoeba} associated to~$\mathfrak a$ to be
\[
\log\vert\mathsf V(\mathfrak a)\vert=\{(\log\vert
z_1\vert,\dots,\log\vert z_d\vert)\mid\mathbf z\in\mathsf
V(\mathfrak a)\}.
\]
Now let~$M$ be a Noetherian~$R_d$ module. Then define
\[
\nesv_k(\alpha_M)=\bigcup_{\mathfrak p\in\ass(M)}\{V\in\mathsf
G_k\mid V^{\perp}\cap\log\vert V(\mathfrak p)\vert\neq\varnothing\}
\]
where~$V^{\perp}$ denotes the orthogonal complement of~$V$
in~$\mathbb R^d$.
The main result in~\cite[Th.~8.4]{MR1869066} says that
\begin{equation*}
\nes_k(\alpha_M)=\nesn_k(\alpha_M)\cup\nesv_k(\alpha_M).
\end{equation*}

\section{Periodic points}\label{section:periodicpoints}

Recall that the dynamical zeta function of a map~$T$ is defined
formally as
\begin{equation}\label{equation:dynamicalzetafunction}
\zeta_T(z)=\exp\sum_{n=1}^{\infty}\frac{z^n}{n}\vert\fix_n(T)\vert.
\end{equation}
If~$\vert\fix_n(T)\vert$ is finite for all~$n\ge1$ and grows at most
exponentially, then~\eqref{equation:dynamicalzetafunction} defines a
complex function in some disc. In our setting there is a
fixed~$\mathbb Z^d$-action~$\alpha$, so write~$\zeta_{\mathbf
n}$ for the zeta function of the map~$\alpha^{\mathbf n}$.
Define~$\rat(\alpha)$ to be the set of~$\mathbf{n}\in\mathbb{Z}^d$
for which~$\zeta_{\mathbf n}$ is a rational function. Notice
that any~$\mathbf n\in\mathbb Z^d$ with the property
that~$\fix_j(\alpha^{\mathbf n})$ is infinite for some~$j\ge1$ is
not a member of~$\rat(\alpha)$.

The simplest non-trivial~$\mathbb Z^2$-action is the
`$\times2,\times3$' system, and the idea behind what follows is
already visible in this example.

\begin{eg}\label{example:firstvisittox2x3example}
Consider the~$\mathbb{Z}^2$-action~$\alpha$ dual to the~$\mathbb
Z^2$-action generated by the commuting maps~$\times 2$ and~$\times
3$ on~$\mathbb{Z}[\frac{1}{6}]$. This is the dynamical system
corresponding to the cyclic~$R_2$-module~$M=R_2/(u_1-2, u_2-3)$. The
set~$\nes_1(\alpha)$ for this example is shown in
Figure~\ref{figurex2x3}; it consists of three lines
with~$\nesn_1(\alpha)$ comprising~$2^{n_1}=1$ and~$3^{n_2}=1$
and~$\nesv_1(\alpha)$ being the single irrational
line~$2^{n_1}3^{n_2}=1$.

\begin{figure}[ht]
\setlength{\unitlength}{0.04cm}
\begin{picture}(100,96)
\thicklines\put(102,48){\tiny$\nesn_1(\alpha)$}
\put(46,96){\tiny$\nesn_1(\alpha)$}\put(102,14){\tiny$\nesv_1(\alpha)$}
\put(50,50){\line(1,0){50}}\put(50,50){\line(3,-2){50}}\put(50,50){\line(-1,0){50}}
\put(50,50){\line(0,1){45}}\put(50,50){\line(0,-1){45}}\put(50,50){\line(-3,2){50}}
\multiput(5,35)(15,0){7}{\circle*{1}}
\multiput(20,20)(15,0){5}{\circle*{1}}
\multiput(20,80)(15,0){5}{\circle*{1}}
\multiput(5,50)(15,0){7}{\circle*{1}}
\multiput(5,65)(15,0){7}{\circle*{1}}
\end{picture}
\caption{\label{figurex2x3}The three non-expansive lines
for~$\times2,\times3$.}
\end{figure}
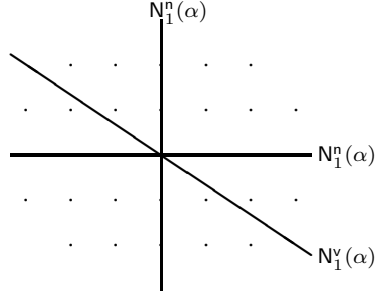

\noindent Figure~\ref{figurex2x3periodicpoints} shows the
map~$\mathbf n\mapsto\vert\fix_{\mathbf n}(\alpha)\vert\in\mathbb N$
for the same system. Notice that~$\vert\fix_{-\mathbf
n}(\alpha)\vert=\vert\fix_{\mathbf n}(\alpha)\vert$, so only the
region~$n_2\ge0$ is shown, with~$\infty$ denoting the
lattice point~$(0,0)$ corresponding to the identity map.

\begin{figure}[ht]
\[
\begin{array}{ccccccccccc}
211&227&235&239&241&121&485&971&1943&3887&7775\\
49&65&73&77&79&5&161&323&647&1295&2591\\
5&11&19&23&25&13&53&107&215&431&863\\
23&7&1&5&7&1&17&35&71&143&287\\
29&13&5&1&1&1&5&11&23&47&95\\
31&5&7&1&1&\infty&1&1&7&5&31\\
\end{array}
\]
\caption{\label{figurex2x3periodicpoints}Periodic point counts
for~$\times2,\times3$.}
\end{figure}

\noindent It may be shown (see~\cite{MR2180241}
and~\cite[Th.~4.7]{miles_zeta}) that
\[
\rat(\alpha)=\{\mathbf n\in\mathbb Z^d\mid n_1n_2\neq0\}
\]
(the issue here is to show that the zeta functions~$\zeta_{(1,0)}$
and~$\zeta_{(0,1)}$ in the two rational non-expansive lines are not
rational). The question addressed in this paper is the following:
does the data in Figure~\ref{figurex2x3periodicpoints} determine the
subdynamical portrait in Figure~\ref{figurex2x3}? As this example
shows, the rationality set~$\rat(\alpha)$ certainly does not
determine~$\nes_1(\alpha)$, so in particular we are asking if the
periodic point data seen along the rational directions can detect
the presence of an irrational non-expansive direction.
\end{eg}

\begin{eg}\label{example:led}
The zero-dimensional analog of the~$\times2,\times3$ system is
Ledrappier's example~\cite{MR512106}, which is the action~$\alpha$
corresponding to the module~$M=R_2/(2,1+u_1+u_2)$. Using the local
structure of~$X_M$ in terms of completions of the function
field~$\mathbb F_2(t)$ and the periodic point formula
from~\cite{MR1461206} we have
\[
\vert\fix_{(n_1,n_2)}(\alpha)\vert=\vert
t^{n_1}(1+t)^{n_2}-1\vert_{\infty}\vert
t^{n_1}(1+t)^{n_2}-1\vert_{t}\vert t^{n_1}(1+t)^{n_2}-1\vert_{1+t},
\]
the three absolute values being given by
\[
\vert r(t)\vert_{t}=2^{-\ord_t(r(t))}, \vert
r(t)\vert_{\infty}=\vert r(t^{-1})\vert_{t}\mbox{ and }\vert
r(t)\vert_{1+t}=2^{-\ord_{1+t}(r(t))}
\]
where~$r(t)\in\mathbb F_2(t)$
(see~\cite{MR1461206},~\cite{MR2031042} or~\cite{MR1882488} for the
details). The set~$\mathsf N_1(\alpha)$ for this example is shown in
Figure~\ref{figurelednex};~$\nesn_1(\alpha)$ comprises the
lines
\[
n_1=0,\medspace n_2=0\mbox{ and }n_1+n_2=0,
\]
while~$\nesv_1(\alpha)$ is
automatically empty since the associated
prime ideal has an empty variety.
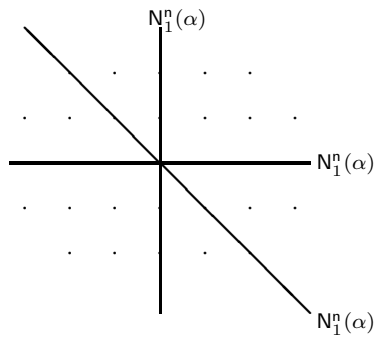
\begin{figure}[ht]
\setlength{\unitlength}{0.04cm}
\begin{picture}(100,96)
\thicklines\put(102,48){\tiny$\nesn_1(\alpha)$}
\put(46,96){\tiny$\nesn_1(\alpha)$}\put(102,-5){\tiny$\nesn_1(\alpha)$}
\put(50,50){\line(1,0){50}}\put(50,50){\line(-1,0){50}}\put(50,50){\line(0,-1){50}}
\put(50,50){\line(0,1){45}}\put(50,50){\line(-1,1){45}}\put(50,50){\line(1,-1){50}}
\multiput(5,35)(15,0){7}{\circle*{1}}
\multiput(20,20)(15,0){5}{\circle*{1}}
\multiput(20,80)(15,0){5}{\circle*{1}}
\multiput(5,50)(15,0){7}{\circle*{1}}
\multiput(5,65)(15,0){7}{\circle*{1}}
\end{picture}
\caption{\label{figurelednex}The three non-expansive lines for
Ledrappier's example.}
\end{figure}

\noindent The periodic point data is shown in
Figure~\ref{figureledperptdata}. Once again the zeta function is
known to be irrational in the non-expansive directions, and in this
example~$\rat(\alpha)$ does indeed detect all the non-expansive
behavior.
\begin{figure}[ht]
\[
\begin{array}{cccccccccccc}
16&32&32&32&32&16&64&128&256&512&1024\\
32&1&16&16&16&1&32&64&128&256&512\\
32&16&4&8&8&4&16&32&64&128&256\\
32&16&8&1&4&1&8&16&32&64&128\\
32&16&8&4&1&1&4&8&16&32&64\\
16&1&4&1&1&\infty&1&1&4&1&16
\end{array}
\]
\caption{\label{figureledperptdata}Periodic point counts for Ledrappier's
example.}
\end{figure}

\noindent In each of the expansive regions, the entries in
Figure~\ref{figureledperptdata} can easily be found using the
ultrametric property of the absolute values. For example, in the
expansive region~$n_1<0$,~$n_2>0$,~$n_1+n_2>0$ we have
\[
\vert t^{n_1}(1+t)^{n_2}-1\vert_{t}=2^{-n_1},\medspace\vert
t^{n_1}(1+t)^{n_2}-1\vert_{\infty}=2^{n_1+n_2}
\]
and
\[
\vert
t^{n_1}(1+t)^{n_2}-1\vert_{t}=1,
\]
giving~$\vert\fix_{(n_1,n_2)}(\alpha)\vert=2^{n_2}.$

The non-expansive lines may be found similarly, though the
resulting expression is a little more involved. For example,
in~\cite[Ex.~8.5]{MR1461206} it is shown
that~$\vert\fix_{(n,0)}(\alpha)\vert=2^{n-2^{\ord_2(n)}}$.

An alternative way to compute the number of periodic points,
better adapted
to more complicated situations, may be found in~\cite[Lem.~4.8]{miles_zeta}.
\end{eg}

To convert the periodic point data into a form which exposes
the expansive subdynamics, we introduce a normalized encoding
of the rational zeta functions arising from elements of the action.

\begin{definition}
Given a rational function~$h\in\mathbb{C}(z)$, denote the set of
poles and zeros of~$h$ by~$\Psi(h)\subset\mathbb{C}$. Let~$\alpha$
be an algebraic~$\mathbb Z^d$-action of entropy rank one, and define
\[
\Omega_\alpha= \left\{ (\hat{\mathbf{n}},|z|^{1/\|\mathbf{n}\|})\mid
z\in\Psi(\zeta_\mathbf{n}), \mathbf{n}\in\rat(\alpha)
\right\}\subset\sph_{d-1}\times\mathbb{R}
\]
where~$\hat{\mathbf{n}}$ denotes the unit vector in the direction
of~$\mathbf{n}$.
\end{definition}

In order to exhibit the relationship between~$\Omega_{\alpha}$
and~$\nes(\alpha)$ we need a `formula' for~$\vert\fix_{\mathbf
n}(\alpha)\vert$, and this has been found by Miles~\cite{miles_zeta}
using the structure of entropy rank one systems
from~\cite{MR2031042}.

A \emph{character} is a homomorphism from an abelian group
into~$\mathbb{C}^\times$. We will be particularly interested in
characters of the form
$\chi:\mathbb{Z}^d\rightarrow\mathbb{C}^\times$. A \emph{real
character} is one with real image. By a \emph{list} we mean a finite
sequence of the form~$L=\langle \chi_1,\dots,\chi_n\rangle$ which
allows for multiplicities. The notation
\[
\chi_L=\chi_1 \chi_2\dots \chi_n
\]
is used to denote the product over all elements of~$L$, with the
understanding that~$\chi_\varnothing\equiv 1$.

Let~$\mathfrak{p}\subset R_d$ be a prime ideal
with~$\textrm{coht}(\mathfrak{p})=1$ and let~$K$ be the field of
fractions of~$R_d/\mathfrak{p}$. Assume
that~$\ch(R_d/\mathfrak{p})=0$, so all the infinite places are
archimedean. These infinite places are uniquely determined by the
embeddings of~$R_d/\mathfrak{p}$ into~$\mathbb{C}$. A point~$z\in
V_\mathbb{C}(\mathfrak{p})$ determines a ring homomorphism
into~$\mathbb{C}$ via the substitution map~$f+\mathfrak{p}\mapsto
f(z)$. The map is injective because~$R_d/\mathfrak{p}$ has Krull
dimension~$1$. Each~$z\in V_\mathbb{C}(\mathfrak{p})$ induces a
character on~$\mathbb{Z}^d$ in an obvious way; there are finitely
many such characters and the coordinates of these are all algebraic
numbers. More generally, any place~$w$ of a domain of the
form~$R_d/\mathfrak{p}$ induces a real character on~$\mathbb{Z}^d$
via the map
\[
(n_1,\dots,n_d)\mapsto(|\overline{u}_i|_w^{n_1},\dots,
|\overline{u}_d|_w^{n_d}),
\]
where~$\overline{u}_i$ denotes the image of~$u_i$
in~$R_d/\mathfrak{p}$,~$i=1\dots d$. This will always be our method
of constructing characters using non-archimedean places.

Using the construction of characters given above, for a
prime~$\mathbb{Z}^d$-action~$\alpha_{R_d/\mathfrak{p}}$
with~$\textrm{coht}(\mathfrak{p})=1$,
let~$\mathcal{W}(R_d/\mathfrak{p})$ be the list of characters
induced by the non-archimedean~$v\in\mathcal{S}_\mathfrak{p}$ and
let~$\mathcal{V}(R_d/\mathfrak{p})$ be the list of characters
induced by the distinct complex embeddings of~$R_d/\mathfrak{p}$.
Note that~$\mathcal{V}(R_d/\mathfrak{p})=\varnothing$
when~$\ch(R_d/\mathfrak{p})>0$.

Now suppose that~$\alpha_M$ is a Noetherian entropy rank one action.
The module~$M$ admits a \emph{prime filtration}
\begin{equation}\label{prime_filtration}
\{0\}=M_0\subset M_1 \subset \dots \subset M_n=M
\end{equation}
where for each~$k$,~$1\leq k\leq n$ we have~$M_k/M_{k-1}\cong
R_d/\mathfrak{q}_k$ for a prime ideal~$\mathfrak{q}_k\subset R_d$
which is either an associated prime of~$M$ or which contains an
associated prime of~$M$. Lemma~8.2 of~\cite{MR2031042} shows that
each minimal element of~$\ass(M)$ always appears in such a
filtration with a fixed multiplicity~$m(\mathfrak{p})$,
so~$m(\mathfrak{p})$ is well-defined for
all~$\mathfrak{p}\in\ass(M)$ with~$\textrm{coht}(\mathfrak{p})=1$.
Set
\[
\mathcal{W}(M)=\bigsqcup
\mathcal{W}(R_d/\mathfrak{p}),
\]
\[
\mathcal{V}(M)=\bigsqcup
\mathcal{V}(R_d/\mathfrak{p}),
\]
where the union of lists is taken over all~$\mathfrak{p}\in\ass(M)$
with~$\textrm{coht}(\mathfrak{p})=1$, ensuring that each
prime~$\mathfrak{p}$ appears with the appropriate
multiplicity~$m(\mathfrak{p})$.

If~$M$ has torsion-free rank one, then~$\mathcal{V}(M)$ has a
particularly simple form.

\begin{lem}\label{torsion_free_rank_one_lemma}
Let~$M$ be a Noetherian~$R_d$-module of torsion-free rank one, and
suppose~$\alpha_M$ has entropy rank one. Then~$\mathcal{V}(M)$
contains one element.
\end{lem}

\begin{proof}
Consider a prime filtration of~$M$ of the
form~\eqref{prime_filtration}. Since~$M$ has torsion-free rank one,
there is at least one associated prime~$\mathfrak{p}$ such
that~$\ch(R_d/\mathfrak{p})=0$. Let~$k\leq n$ be the least integer
such that~$\ch(R_d/\mathfrak{q}_k)=0$.
Then~$\coht(\mathfrak{q}_k)=1$ and~$\mathfrak{q}_k\in\ass(M)$.
Suppose~$k<n$. Since~$M_{k+1}/M_k\cong R_d/\mathfrak{q}_{k+1}$,
there exists~$a\in M_{k+1}\setminus M_k$ such that any element
of~$M_{k+1}$ can be written in the form~$x+fa$ for some~$x\in M_k$
and~$f\in R_d$ with~$fa\in M_k$ if and only
if~$f\in\mathfrak{q}_{k+1}$. However, both~$M_k$ and~$M$ have
torsion-free rank one so there exists~$c\in\mathbb{Z}$ such
that~$ca\in M_k$. Therefore,~$c\in\mathfrak{q}_{k+1}$
and~$\ch(R_d/\mathfrak{q}_{k+1})>0$. In a similar way, it follows
that~$\ch(R_d/\mathfrak{q}_{j})>0$ for all~$j>k$.
Hence~$\mathfrak{q}_k$ is the only prime
with~$\ch(R_d/\mathfrak{q}_k)=0$; moreover~$m(\mathfrak{q}_k)=1$.
So
\[
\mathcal{V}(M)=\mathcal{V}(R_d/\mathfrak{q}_k).
\]
If~$k=n$ then
again~$\mathcal{V}(M)=\mathcal{V}(R_d/\mathfrak{q}_k)$.
Finally,~$R_d/\mathfrak{q}_k$ is isomorphic to a subring
of~$\mathbb{Q}$, so~$\mathcal{V}(M)$ contains one character induced
by the single infinite place of~$\mathbb{Q}$.
\end{proof}

Any character~$\chi:\mathbb{Z}^d\rightarrow\mathbb{C}^\times$
induces a real character~$\chi^*$ on~$\mathbb{R}^d$, by setting
\[
\chi^*(\kappa\mathbf{e}_i)=|\chi(\mathbf{e}_i)|^{\kappa}
\]
where~$\kappa\in\mathbb{R}$ and~$\mathbf{e}_i$ is the
standard~$i$-th basis vector in~$\mathbb{Z}^d$,~$i=1\dots d$.
Applying this construction to an element of~$\mathcal{W}(M)$ yields
a genuine extension, but the same is not necessarily true for
elements of~$\mathcal{V}(M)$.

\begin{prop}\label{non_expansive_hyperplanes_description}
Suppose~$\alpha_M$ is an algebraic~$\mathbb{Z}^d$-action of
expansive rank one. Then~$\nes_{d-1}(\alpha_M)$ consists precisely
of the finite set of hyperplanes defined by the equations
\[
\chi^*(\mathbf{n})=1
\]
where~$\chi\in\mathcal{V}(M)\cup\mathcal{W}(M)$.
Furthermore,~$\nesv_{d-1}(\alpha_M)$ is determined by those
characters in~$\mathcal{V}(M)$ and~$\nesn_{d-1}(\alpha_M)$ by those
characters in~$\mathcal{W}(M)$.
\end{prop}

\begin{proof}
This is a combination of~\cite[Th.~4.3.10]{miles_thesis}
and~\cite[Th.~8.4]{MR1869066}.
\end{proof}

By expressing~$\nes_1(\alpha_M)$ in terms of the intersection of
non-expansive lines with~$\sph_{d-1}$ and  referring
to~\cite[Th.~3.6]{MR1355295}, we also find the following description
of the expansive subdynamics.

\begin{cor}\label{non_expansive_lines_description}
If~$\alpha_M$ is an algebraic~$\mathbb{Z}^d$-action of expansive
rank one then
\[
\nesv_1(\alpha) = \bigcup_{\chi\in\mathcal{V}(M)}\{\mathbf{v}\in
\sph_{d-1}\mid\chi^*(\mathbf{v})=1\},
\]
\[
\nesn_1(\alpha) = \bigcup_{\chi\in\mathcal{W}(M)}\{\mathbf{v}\in
\sph_{d-1}\mid\chi^*(\mathbf{v})=1\}
\]
and the set of expansive directions is dense in~$\sph_{d-1}$.
\end{cor}

\section{Main Results}

To begin this section, we return to
the examples in Section~\ref{section:periodicpoints}

\begin{eg}\label{2_3_example}
Let~$\alpha$ be the~$\mathbb Z^2$-action corresponding to the
$R_2$-module
\[
M=R_2/(u_1-2, u_2-3)
\]
discussed in Example~\ref{example:firstvisittox2x3example}.
Recall that
\[
\rat(\alpha)=\{\mathbf n\in\mathbb Z^d\mid n_1n_2\neq0\}.
\]
Here,~$\rat(\alpha)$ consists precisely of
those~$\mathbf{n}\in\mathbb{Z}^2$ for which~$\alpha_M^{\mathbf{n}}$
is expansive, but this need not always be the case
(see~\cite[Ex.~4.3]{miles_zeta} for an example). Notice that
expansiveness of the elements~$\alpha_M^{\mathbf{n}}$ of the action
can only ever detect non-expansiveness in rational directions, so
the irrational line in~$\nes_1(\alpha)$ will be missed.
For~$\mathbf{n}=(n_1,n_2)\in\rat(\alpha)$, using the periodic point
formula from~\cite{miles_zeta},
\[
|\fix_j(\alpha_M^{\mathbf{n}})|= |2^{jn_1}3^{jn_2}-1|_\infty
|2^{jn_1}3^{jn_2}-1|_2 |2^{jn_1}3^{jn_2}-1|_3.
\]
It follows that
\[
\zeta_\mathbf{n}(z)=(1-g(\mathbf{n})z)^{\lambda_1}(1-g(\mathbf{n})2^{n_1}3^{n_2}z)^{\lambda_1},
\]
where~$\lambda_1,\lambda_2\in\{-1,1\}$ and
\[
g(\mathbf{n})=|2^{n_1}-1|_2|3^{n_2}-1|_3.
\]

\begin{figure}[!h]
\centering
\psfrag{ x0}{$0$}
\psfrag{ x1}{$1$}
\psfrag{ x2}{$2$}
\psfrag{ x3}{$3$}
\psfrag{ x4}{$4$}
\psfrag{ x5}{$5$}
\psfrag{ x6}{$6$}
\psfrag{ y1}{$1$}
\psfrag{ y0.8}{$0.8$}
\psfrag{ y0.6}{$0.6$}
\psfrag{ y0.4}{$0.4$}
\psfrag{ y0.2}{$0.2$}
\psfrag{ y0}{$0$}
\includegraphics[scale=0.75]{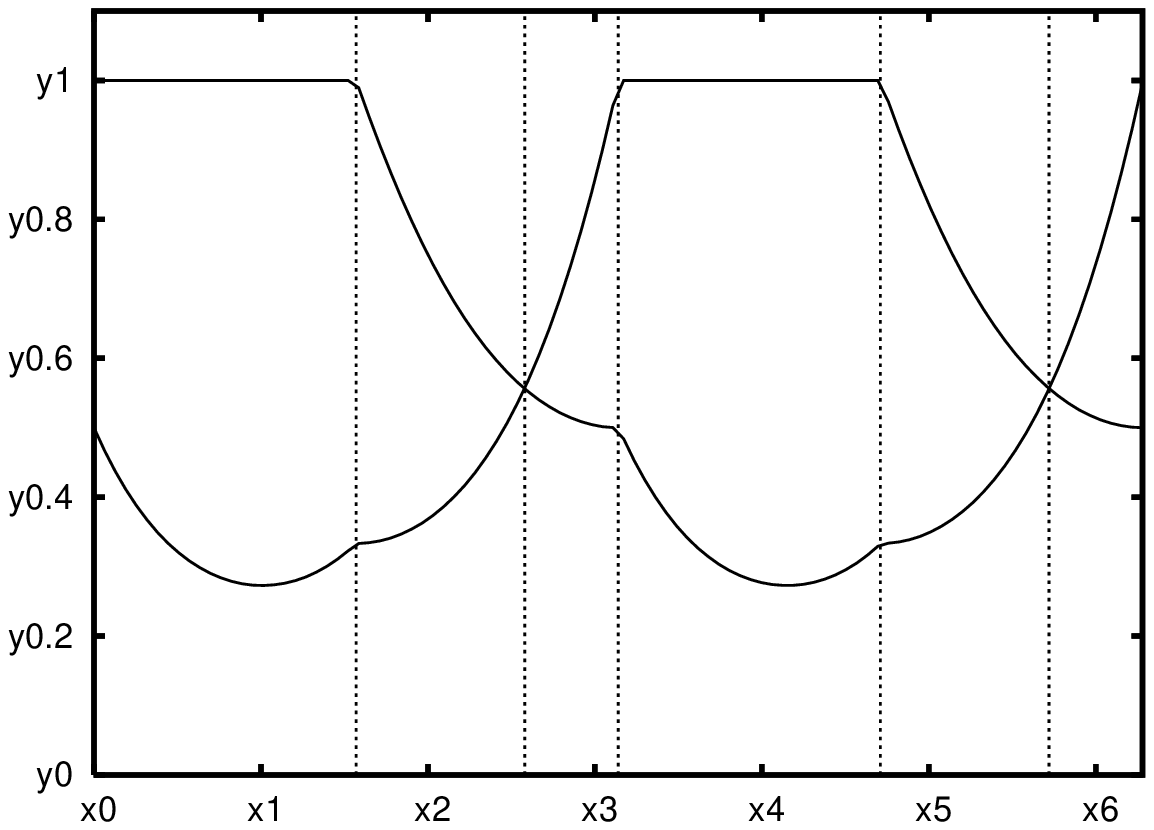}
\caption{$\overline{\Omega}_\alpha$ for the~$\times2,\times3$ system in Example~\ref{2_3_example}.}
\label{2_3_figure}
\end{figure}

\noindent The resulting directional pole and zero
data~$\overline{\Omega}_\alpha$, realized as a subset
of~$[0,2\pi)\times\mathbb{R}$, is shown in Figure~\ref{2_3_figure}.
Non-expansive directions are marked with a dashed line.
\end{eg}

\begin{eg}\label{example:led2}
Let~$\alpha$ be the~$\mathbb Z^2$-action corresponding to the
$R_2$-module
\[
M=R_2/(2,1+u_1+u_2)
\]
discussed in Example~\ref{example:led}
Recall that
\[
\rat(\alpha)=\{\mathbf n\in\mathbb Z^d\mid n_1n_2\neq0
\mbox{ and }n_1+n_2\neq 0\}.
\]
For~$\mathbf{n}=(n_1,n_2)\in\rat(\alpha)$, using the periodic point
formula from~\cite{miles_zeta},
\[
\zeta_\mathbf{n}(z)=
(1-g(\mathbf{n})z)^{-1}
\]
where
\[
g(\mathbf{n})=|t^{n_1}(1+t)^{n_2}-1|_{1+t}\vert t^{n_1}(1+t)^{n_2}-1|_{t}.
\]
\begin{figure}[!h]
\centering
\psfrag{0.3}{$0.3$}
\psfrag{0.35}{}
\psfrag{0.4}{$0.4$}
\psfrag{0.45}{}
\psfrag{0.5}{$0.5$}
\psfrag{0.55}{}
\psfrag{0.6}{$0.6$}
\psfrag{0.65}{}
\psfrag{ 0}{$0$}
\psfrag{ 1}{$1$}
\psfrag{ 2}{$2$}
\psfrag{ 3}{$3$}
\psfrag{ 4}{$4$}
\psfrag{ 5}{$5$}
\psfrag{ 6}{$6$}
\includegraphics[scale=0.75]{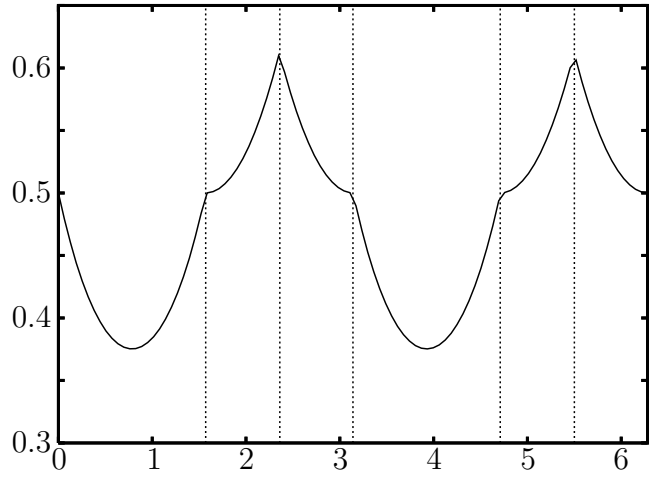}
\caption{$\overline{\Omega}_\alpha$ for Ledrappier's example.}
\label{ledzetafigure}
\end{figure}
The resulting directional pole and zero
data~$\overline{\Omega}_\alpha$, realized as a subset
of~$[0,2\pi)\times\mathbb{R}$, is shown in Figure~\ref{ledzetafigure}.
Non-expansive directions are marked with a dashed line.
\end{eg}

In Examples~\ref{2_3_example} and~\ref{example:led2}
the periodic point `formula' is
particularly simple, mainly because the module~$M$ is cyclic. In
order to deal with more general modules, some machinery is needed.

\begin{prop}\label{periodic_point_structure_prop}
Let~$\alpha_M$ be a Noetherian
algebraic~$\mathbb{Z}^d$-action,~$j\in\mathbb{N}$
and~$\mathbf{n}\in\mathbb{Z}^d$. The following conditions are
equivalent.
\begin{enumerate}
\item $\alpha_M^\mathbf{n}$ has finitely many points of period~$j$.
\item $\alpha_{R_d/\mathfrak{p}}^\mathbf{n}$ has finitely many points of period~$j$ for
every~$\mathfrak{p}\in\ass(M)$.
\end{enumerate}
\end{prop}

For the proof, the following lemma is needed.

\begin{lem}\label{periodic_point_counting_lemma}
Let~$j\in\mathbb{N}$,~$\mathbf{n}\in\mathbb{Z}^d$,~$\mathfrak{b}=(u^{j\mathbf{n}}-1)\subset
R_d$ and let~$N$ be an~$R_d$-module.
\begin{enumerate}
\item Whenever~$N/\mathfrak{b}N$ or~$|F_j(\alpha^{\mathbf{n}}_N)|$ is finite,
\[
\left|F_j(\alpha^{\mathbf{n}}_N)\right| = |N/\mathfrak{b}N|.
\]
\item If~$L\subset N$ is a submodule satisfying~$N/L\cong R_d/\mathfrak{q}$
for some prime ideal~$\mathfrak{q}\subset R_d$, then
\[
\mathfrak{b}N\cap L =
\left\{
\begin{array}{ll}
\mathfrak{b}N & \textrm{when }u^{j\mathbf{n}}-1\in\mathfrak{q},\\
\mathfrak{b}L & \textrm{otherwise}.
\end{array}
\right.
\]
\end{enumerate}
\end{lem}

\begin{proof}
Part~(1) is a simpler version of (for
example)~\cite[Lem.~7.2]{MR1062797}; the hypothesis means
that~$F_j(\alpha^{\mathbf{n}}_N)$ is a finite abelian group and the
statement just says that such a group has the same number of
elements as its character group, which is well known. Part~(2) is a
modification of~\cite[Lem.~3.4]{miles_zeta} which uses an ergodicity
assumption to get a slightly stronger result, but this is not needed
here.
\end{proof}

\begin{proof}[Proof of Proposition \ref{periodic_point_structure_prop}]
Set~$\mathfrak{b}=(u^{j\mathbf{n}}-1)$ and let~$L\subset N$ be
submodules of~$M$ satisfying~$N/L\cong R_d/\mathfrak{q}$ for a prime
ideal~$\mathfrak{q}\subset R_d$. Then by
Lemma~\ref{periodic_point_counting_lemma}(2),
\begin{eqnarray}
|N/\mathfrak{b}N|
 & = &
\left|\frac{N/L}{(\mathfrak{b}N+L)/L}\right|
\left|\frac{\mathfrak{b}N+L}{\mathfrak{b}N}\right|\nonumber \\
 & = &
\left|\frac{N/L}{\mathfrak{b}(N/L)}\right|
\left|\frac{L}{\mathfrak{b}N\cap L}\right|\nonumber \\
 & = &
\left\{
\begin{array}{ll}
|N/L|
\left|\displaystyle\frac{L/\mathfrak{b}L}{\mathfrak{b}N/\mathfrak{b}L}\right|
 & \textrm{when }\mathfrak{b}\subset\mathfrak{q}, \\
\left|\displaystyle\frac{N/L}{\mathfrak{b}(N/L)}\right|
\left|\displaystyle\frac{L}{\mathfrak{b}L}\right|
 & \textrm{otherwise}.
\end{array}
\right. \label{module_decomposition_equation}
\end{eqnarray}
To prove the result, we will first show
\begin{equation}\label{induction_step_equation}
|N/\mathfrak{b}N|<\infty
\Leftrightarrow
|L/\mathfrak{b}L|<\infty
\textrm{ and }
\left|\frac{N/L}{\mathfrak{b}(N/L)}\right|<\infty,
\end{equation}
and then proceed by induction.

($\Rightarrow$) If~$\mathfrak{b}\not\subset\mathfrak{q}$ then the
required implication follows immediately from
equation~\eqref{module_decomposition_equation}. Assume now
that~$\mathfrak{b}\subset\mathfrak{q}$. The finiteness of
equation~\eqref{module_decomposition_equation} forces~$N/L$ to be
finite and hence~$|(N/L)/\mathfrak{b}(N/L)|<\infty$. To see
that~$|L/\mathfrak{b}L|<\infty$,
let~$\phi:N\rightarrow~\mathfrak{b}N/\mathfrak{b}L$ be defined by
multiplication by~$u^{j\mathbf{n}}-1$, followed by the natural
quotient map. Clearly~$\phi$ is surjective
and~$L\subset\textrm{ker}(\phi)$ which means there is a surjective
homomorphism from~$N/L$ to~$\mathfrak{b}N/\mathfrak{b}L$.
Hence~$\mathfrak{b}N/\mathfrak{b}L$ is finite and~$L/\mathfrak{b}L$
must also be finite to ensure that
equation~\eqref{module_decomposition_equation} holds.

($\Leftarrow$) If~$\mathfrak{b}\not\subset\mathfrak{q}$ then
equation~\eqref{module_decomposition_equation} immediately
gives~$|N/\mathfrak{b}L|<\infty$. For the
case~$\mathfrak{b}\subset\mathfrak{q}$, again apply
equation~\eqref{module_decomposition_equation}, noting that
\[
|N/L| = \left|\frac{N/L}{\mathfrak{b}(N/L)}\right| \mbox{ and }
\left|\frac{L/\mathfrak{b}L}{\mathfrak{b}N/\mathfrak{b}L}\right|
\leq |L/\mathfrak{b}L|.
\]

($1\Rightarrow 2$) A prime filtration of~$M$ of the
form~\eqref{prime_filtration} may be taken so that any
chosen~$\mathfrak{p}\in\ass(M)$ appears as~$\mathfrak{q}_1$ (see for
example the proof of~\cite[Prop.~3.7]{MR1322960}). Starting
with~$M_n=M$ and descending through the submodules appearing
in~\eqref{prime_filtration}, successively
applying~\eqref{induction_step_equation}
gives~$|M_1/\mathfrak{b}M_1|<\infty$. The required result follows by
Lemma~\ref{periodic_point_counting_lemma}(1).

($2\Rightarrow 1$) Again take a prime filtration of~$M$ of the
form~\eqref{prime_filtration}. For each~$1\leq k\leq
n$,~$\mathfrak{q}_k$ is an associated prime of~$M$ or contains such
a prime. Let~$\mathfrak{p}\in\ass(M)$  and suppose
that~$\mathfrak{q}_k\supset\mathfrak{p}$. There is a natural
surjective homomorphism
\[
\frac{R_d/\mathfrak{p}}{\mathfrak{b}(R_d/\mathfrak{p})}
\cong
\frac{R_d}{\mathfrak{b}+\mathfrak{p}}
\longrightarrow
\frac{R_d}{\mathfrak{b}+\mathfrak{q}_k}
\cong
\frac{R_d/\mathfrak{q}_k}{\mathfrak{b}(R_d/\mathfrak{q}_k)}.
\]
Since~$|F_j(\alpha^{\mathbf{n}}_{R_d/\mathfrak{p}})|<\infty$ for
all~$\mathfrak{p}\in\ass(M)$, by the above and part~(1) of
Lemma~\ref{periodic_point_counting_lemma},
\[
\left|\frac{M_k/M_{k-1}}{\mathfrak{b}(M_k/M_{k-1})}\right|<\infty
\]
for all~$1\leq k\leq n$. Starting with~$k=1$ and ascending through
the modules appearing in~\eqref{prime_filtration}, successively
applying~\eqref{induction_step_equation}, gives
\[
|M_n/\mathfrak{b}M_n|<\infty
\]
and hence~$|F_j(\alpha^{\mathbf{n}}_M)|<\infty$.
\end{proof}

The following consequence of
Proposition~\ref{periodic_point_structure_prop} shows that for
Noetherian actions with entropy rank greater than one, periodic
point data for individual elements is much less useful than in the
rank one case.

\begin{cor}
Let~$\alpha_M$ be a Noetherian algebraic~$\mathbb{Z}^d$-action with
entropy rank greater than one. For each~$\mathbf{n}\in\mathbb{Z}^d$,
the automorphism~$\alpha_M^\mathbf{n}$ has infinitely many points of
any given period.
\end{cor}

\begin{proof}
First note, by assumption~$d>1$. Since~$\entrk(\alpha_M)>1$, by
Proposition~\ref{main_structure_theorem}, there is at least
one~$\mathfrak p\in\ass(M)$ with~$\coht(\mathfrak{p})>1$.
Let~$j\in\mathbb{N}$,~$D=R_d/\mathfrak p$ and
consider~$|F_j(\alpha_{R_d/\mathfrak p}^\mathbf{n})|$. Suppose this
is finite. Then by Lemma~\ref{periodic_point_counting_lemma}(1),
there is a proper principal ideal~$\mathfrak a \subset D$ such that
\[
|F_j(\alpha_{R_d/\mathfrak p}^\mathbf{n})|=|D/\mathfrak a|.
\]
By the Principal Ideal Theorem~\cite[Th.~10.1]{MR1322960} there is a
prime ideal~$\mathfrak{q}\subset D$ containing~$\mathfrak a$
with~$\hght(\mathfrak q)=1$. Since~$D$ is a quotient of~$R_d$ by a
prime ideal of coheight at least~$2$ and since all maximal ideals
of~$R_d$ have the same height, every maximal ideal of~$D$ has height
at least~$2$. Hence~$\kdim(D/\mathfrak q)=\coht(\mathfrak q)\geq 1$.
Therefore,~$D/\mathfrak q$ cannot be finite as it either
contains~$\mathbb{Z}$ or has non-zero transcendence degree over a
finite field. Because there is a surjection from~$D/\mathfrak a$
to~$D/\mathfrak q$, the supposition that~$|F_j(\alpha_{R_d/\mathfrak
p}^\mathbf{n})|<\infty$ must therefore be false. The result now
follows by Proposition~\ref{periodic_point_structure_prop}.
\end{proof}

Recall that an algebraic~$\mathbb Z^d$-action~$\alpha$ is ergodic
(always with respect to Haar measure) if there are no non-trivial
invariant sets for the whole action;
Schmidt~\cite[Th.~6.5]{MR1345152} shows that~$\alpha$ is ergodic if
and only if there is some~$\mathbf n\in\mathbb Z^d$ for which the
single automorphism~$\alpha^{\mathbf n}$ is ergodic. The next result
gives a criteria to identify certain elements of the action that
must be ergodic. Notice that an automorphism of a finite group
cannot be ergodic since it must fix the identity, which has positive
Haar measure.

\begin{cor}\label{periodic_points_and_ergodicity_result}
Suppose~$\alpha_M$ is an ergodic Noetherian
algebraic~$\mathbb{Z}^d$-action and~$\mathbf{n}\in\mathbb{Z}^d$. If
the automorphism~$\alpha_M^\mathbf{n}$ has finitely many points of
every period then~$\alpha_M^\mathbf{n}$ is ergodic.
\end{cor}

\begin{proof}
Since~$\alpha_M$ is ergodic,~\cite[Th.~6.5]{MR1345152} shows that
there exists~$\mathbf{m}\in\mathbb{Z}^d$ such
that~$\alpha_{R_d/\mathfrak{p}}^\mathbf{m}$ is ergodic for
every~$\mathfrak{p}\in\ass(M)$. If~$\mathfrak{p}$ is maximal
then~$\alpha^\mathbf{m}_{R_d/\mathfrak{p}}$ is an automorphism of a
finite field which cannot be ergodic. Hence~$\ass(M)$ contains no
maximal ideals.

For a contradiction suppose that~$\alpha_M^{\mathbf{n}}$ is not
ergodic. Then, again using~\cite[Th.~6.5]{MR1345152}, there
exists~$j\geq 1$ and a prime ideal~$\mathfrak{p}\in\ass(M)$ such
that~$u^{j\mathbf{n}}-1\in\mathfrak{p}$. Therefore, the
ideal~$(u^{j\mathbf{n}}-1)\subset R_d$
annihilates~$R_d/\mathfrak{p}$ and by
Lemma~\ref{periodic_point_counting_lemma}(1),
\[
|F_j(\alpha_{R_d/\mathfrak{p}}^\mathbf{n})|=|R_d/\mathfrak{p}|.
\]
Since~$\mathfrak{p}$ is not maximal this is not finite, which
contradicts Proposition~\ref{periodic_point_structure_prop}.
\end{proof}

This allows us to describe the set~$\Psi(\zeta_{\mathbf n})$ in
terms of places and characters.

\begin{prop}\label{coefficient_descriptions}
Suppose~$\alpha_M$ is an ergodic Noetherian
algebraic~$\mathbb{Z}^d$-action of entropy rank one.
If~$\mathbf{n}\in\rat(\alpha_M)$ then
\begin{equation}\label{pole_and_zero_set}
\Psi(\zeta_\mathbf{n})
 =
\left\{
\mu_{\mathbf{n}}\chi_T\chi_L(-\mathbf{n})\mid L\subset\mathcal{V}(M)
\right\}
\end{equation}
where~$T=\langle\chi\in\mathcal{W}(M)\mid\chi(\mathbf{n})>1\rangle$
and~$\mu_{\mathbf{n}}\in\{-1,1\}$.
\end{prop}

\begin{proof}
Since~$\mathbf{n}\in\rat(\alpha_M)$,~$|F_j(\alpha_M^\mathbf{n})|<\infty$
for all~$j\in\mathbb{N}$. Therefore, by
Corollary~\ref{periodic_points_and_ergodicity_result}, we may assume
that~$\alpha_M^\mathbf{n}$ is ergodic. Hence, the explicit formula
for~$\zeta_\mathbf{n}$ given in~\cite{miles_zeta} may be used. Let
\[
R=\langle \chi\in\mathcal{V}(M)\mid |\chi(\mathbf{n})|>1\rangle
\]
and denote the complement of~$J\subset\mathcal{V}(M)$
in~$\mathcal{V}(M)$ by~$J'$. Applying~\cite[Th.~4.7,
Lem.~4.2]{miles_zeta}, if~$\zeta_{\mathbf{n}}$ is rational then
\[
\zeta_\mathbf{n}(z)
=
\prod_{J\subset\mathcal{V}(M)}(1-c_Jz)^{\lambda_J},
\]
where
\begin{equation}\label{zeta_coefficient_from_miles}
c_J=\mu_{\mathbf{n}}\chi_T\chi_{R\cap J'}\chi_{R'\cap J}(\mathbf{n})
\end{equation}
and~$\mu_{\mathbf{n}},\lambda_J\in\{-1,1\}$. Hence,
setting~$L=(R\cap J')\cup(R'\cap J)$ shows that the corresponding
pole or zero of~$\zeta_\mathbf{n}$ arising from the
coefficient~\eqref{zeta_coefficient_from_miles} is accounted for
in~\eqref{pole_and_zero_set}. Conversely, to see that for a
given~$L\subset\mathcal{V}(M)$,~$\mu_{\mathbf{n}}\chi_T\chi_L(-\mathbf{n})\in\Psi(\zeta_\mathbf{n})$,
simply set~$J=(R\cap L')\cup(R'\cap L)$
in~\eqref{zeta_coefficient_from_miles}.
\end{proof}

We can now prove the main result, which roughly says
that~$\overline{\Omega}_\alpha$ (which is determined completely by
the periodic point data for~$\alpha$) is in fact a set of graphs, and that
these graphs pick out two sets of non-expansive directions. The first is
the \emph{crossing set}, where graphs of functions cross, and this
turns out to contain~$\nesv_1(\alpha)$. The second is the set of
points where one of the functions is not differentiable, and this
turns out to be exactly~$\nesn_1(\alpha)$. In the simplest
situations the crossing set is exactly~$\nesv_1(\alpha)$, and the
extra directions identified in the general case reflect vanishing of
other combinations of Lyapunov exponents.

\begin{thm}\label{main_theorem}
Let~$\alpha$ be an ergodic algebraic~$\mathbb{Z}^d$-action of
expansive rank one. Then there exist continuous
functions~$f_1,\dots,f_r\in\mathcal{C}(\sph_{d-1},\mathbb{R})$
such that
\[
\overline{\Omega}_\alpha
 =
\bigcup_{k=1}^{r}
\left\{(\mathbf{v},f_k(\mathbf{v}))\mid\mathbf{v}\in\sph_{d-1}\right\}.
\]
Moreover,~$f_1,\dots,f_r$ have the property that
\begin{enumerate}
\item\label{variety_non_exp_relation}
$\nesv_1(\alpha) \subset \{\mathbf{v}\in\sph_{d-1}\mid
f_j(\mathbf{v})=f_k(\mathbf{v}) \mbox{ for some }j\neq k\}$,
\item\label{noetherian_non_exp_relation}
$\nesn_1(\alpha) = \{\mathbf{v}\in\sph_{d-1}\mid f_k \mbox{ is
not smooth at } \mathbf{v} \mbox{ for some }k\}$.
\end{enumerate}
If~$\dim(X)\leq 1$ then equality holds
in~\eqref{variety_non_exp_relation}.
\end{thm}

\begin{proof}
Let~$\alpha=\alpha_M$. First note that
given~$\chi\in\mathcal{V}(M)\cup\mathcal{W}(M)$, the induced
character~$\chi^*$ is continuous on~$\mathbb{R}^d$. For
each~$\chi\in\mathcal{W}(M)$,
define~$g_\chi:\mathbb{R}^d\rightarrow\mathbb{R}^\times$ by
\[
g_\chi(\mathbf{n})=\max\{\chi^*(\mathbf{n}),1\}.
\]
Then~$g=\prod_{\chi\in\mathcal{W}(M)}g_\chi$ is continuous
since~$\mathcal W(M)$ is finite. Now, for
each~$L\subset\mathcal{V}(M)$
define~$f_L:\mathbb{R}^d\rightarrow\mathbb{R}^\times$ by
\[
f_L(\mathbf{n})=g(-\mathbf{n})\chi^*_L(-\mathbf{n}).
\]
Let~$\mathbf{n}\in\mathbb{R}^d$ and define
\[
T(\mathbf{n})=\{\chi\in\mathcal{W}(M)\mid\chi^*(\mathbf{n})>1\}
\]
Then~$T(\mathbf{n})=T(\hat{\mathbf{n}})$ and
\[
\chi^*_{T(\hat{\mathbf{n}})}\chi^*_L(-\hat{\mathbf{n}})=f_L(\hat{\mathbf{n}}).
\]
Thus, if~$\mathbf{n}\in\rat(\alpha)$ then
Proposition~\ref{coefficient_descriptions} shows that
\begin{equation}\label{poles_and_zeros_for_direction_equation}
\left\{
(\hat{\mathbf{n}},|z|^{1/\|\mathbf{n}\|})\mid
z\in\Psi(\zeta_\mathbf{n})
\right\}
 =
\left\{
(\hat{\mathbf{n}},f_L(\hat{\mathbf{n}}))\mid
L\subset\mathcal{V}(M)
\right\}.
\end{equation}
Since~$\alpha$ has expansive rank one, by
Corollary~\ref{non_expansive_lines_description}, the intersection of
the set of expansive lines with~$\sph_{d-1}$ is dense.
Since~$\alpha$ is ergodic,~\cite[Th.~4.4]{miles_zeta} shows that
each expansive automorphism has a rational zeta function. Hence, the
projection of~$\rat(\alpha)$ to~$\sph_{d-1}$ is also dense.
Using~\eqref{poles_and_zeros_for_direction_equation}, it follows
that
\[
\overline{\Omega}_\alpha
=
\bigcup_{L\subset\mathcal{V}(M)}
\left\{
(\mathbf{v},f_L(\mathbf{v}))\mid\mathbf{v}\in\sph_{d-1}
\right\}.
\]
After making the obvious restriction of~$f_L$ to unit vectors,
for~$L\subset\mathcal{V}(M)$, the first part of the theorem follows.

To prove~\eqref{variety_non_exp_relation}, first note by
Corollary~\ref{non_expansive_lines_description}
\begin{equation}\label{chars_def_variety_non_exp_dirs_relation}
\nesv_1(\alpha_M) = \bigcup_{\chi\in\mathcal{V}(M)}\{\mathbf{v}\in\sph_{d-1}
\mid
\chi^*(\mathbf{v})=1\}.
\end{equation}
Furthermore, since~$g(\mathbf{v})\neq 0$ for
all~$\mathbf{v}\in\sph_{d-1}$, for any~$\chi\in
\mathcal{V}(M)$,
\begin{equation}\label{chars_and_f_equal_one_equation}
\chi^*(\mathbf{v})=1
\Leftrightarrow
f_{\{\chi\}}(\mathbf{v})=f_\varnothing(\mathbf{v}).
\end{equation}

To establish~\eqref{noetherian_non_exp_relation}, again by
Corollary~\ref{non_expansive_lines_description}
\[
\nesn_1(\alpha_M) = \bigcup_{\chi\in\mathcal{W}(M)}\{\mathbf{v}\in\sph_{d-1}
\mid\chi^*(\mathbf{v})=1\}.
\]
Since for each~$\chi\in\mathcal{V}(M)$,~$\chi^*$ is smooth,~$f_L$
fails to be so only at points where~$g$ is not smooth. It follows
from the definition of~$g$ that this happens precisely
when~$\chi^*(\mathbf{v})=1$ for some~$\chi\in\mathcal{W}(M)$
and~$\mathbf{v}\in\sph_{d-1}$.

Finally, we show that equality holds
in~\eqref{variety_non_exp_relation} if~$\dim(X)\leq 1$.

If~$\dim(X)=0$ then~$\ch(R_d/\mathfrak{p})>0$ for
all~$\mathfrak{p}\in\ass(M)$. Hence, both~$\mathcal{V}(M)$
and~$\nesv_1(\alpha)$ are empty and equality
in~\eqref{variety_non_exp_relation} holds trivially.

If~$\dim(X)=1$ then Lemma~\ref{torsion_free_rank_one_lemma} shows
that~$\mathcal{V}(M)$ consists of a single character~$\chi$. Hence
for any~$\mathbf{v}\in\sph_{d-1}$, the only possible relation
of the form~$f_J(\mathbf{v})=f_L(\mathbf{v})$
for~$J,L\subset\mathcal{V}(M)$,~$J\neq L$, is given
by~\eqref{chars_and_f_equal_one_equation}. Therefore, equality
in~\eqref{variety_non_exp_relation} follows
from~\eqref{chars_def_variety_non_exp_dirs_relation}.
\end{proof}

\begin{rmk}
The set~$\overline{\Omega}_\alpha$ carries information about the
growth rate of periodic points in all the expansive directions, and
therefore knowledge of the directional entropies in the sense of
Einsiedler and Lind~\cite[Prop.~8.5]{MR2031042}. In the notation of
Theorem~\ref{main_theorem},
\[
h(\alpha^{\mathbf n})=\|\mathbf n\| \log \min_{1\le i\le
r}\{f_i(\hat{\mathbf n})\}.
\]
\end{rmk}

\section{Further examples}

The following example illustrates that equality cannot hold in
Theorem~\ref{main_theorem}\eqref{variety_non_exp_relation} in
general.

\begin{eg}\label{(1+sqrt(2))_(2+sqrt(3))_example}
The following example is algebraically conjugate to
the~$\mathbb{Z}^2$-action of the~$4$-torus
in~\cite[Ex.~6.6]{MR2031042}. Let~$\alpha=\alpha_M$ be
the~$\mathbb{Z}^2$-action corresponding to the
prime~$R_2$-module~$R_2/\mathfrak{p}$, where
\[
\mathfrak{p}=(u_1^2-2u_1-1,u_2^2-4u_2+1)\subset R_2.
\]
The zeros of the first polynomial are~$1\pm\sqrt{2}$ and of the
second are~$2\pm\sqrt{3}$.
Hence~$M=R_2/\mathfrak{p}\cong\mathbb{Z}[\sqrt{2},\sqrt{3}]$
and~$\mathcal{S}_\mathfrak{p}$ consists of~$4$ infinite places
induced by the elements
of~$G=\textrm{Gal}(\mathbb{Q}(\sqrt{2},\sqrt{3})|\mathbb{Q})$.
Therefore
\[
\mathcal{V}(M)=\{(\sigma(1+\sqrt{2}),\sigma(2+\sqrt{3}))\mid\sigma\in G\}
\]
and~$\mathcal{V}(M)$ may be used to establish the poles and zeros
of~$\zeta_\mathbf{n}$ for all~$\mathbf{n}$ in~$\rat(\alpha)$. Note
in this example~$\mathcal{W}(M)=\varnothing$. Consider
\[
L=\{(1+\sqrt{2},2+\sqrt{3}), (1+\sqrt{2},2-\sqrt{3})\}
\subset
\mathcal{V}(M).
\]

The equation
\[
\chi^*_L(\mathbf{v})=\chi^*_\varnothing(\mathbf{v})=1
\]
has a solution~$\mathbf{v}=(0,1)$. However, for
all~$\chi\in\mathcal{V}(M)$,~$\chi((0,1))\neq 1$ and
clearly~$(0,1)\not\in\nes_1(\alpha)$.
\begin{figure}[!h]
\centering
\psfrag{y0}{$0$}
\psfrag{y1}{$1$}
\psfrag{y2}{$2$}
\psfrag{y3}{$3$}
\psfrag{y4}{$4$}
\psfrag{y5}{$5$}
\psfrag{y6}{$6$}
\psfrag{x0}{$0$}
\psfrag{x2}{$2$}
\psfrag{x4}{$4$}
\psfrag{x6}{$6$}
\psfrag{x8}{$8$}
\psfrag{x10}{$10$}
\psfrag{x12}{$12$}
\psfrag{x14}{$14$}
\includegraphics[scale=0.75]{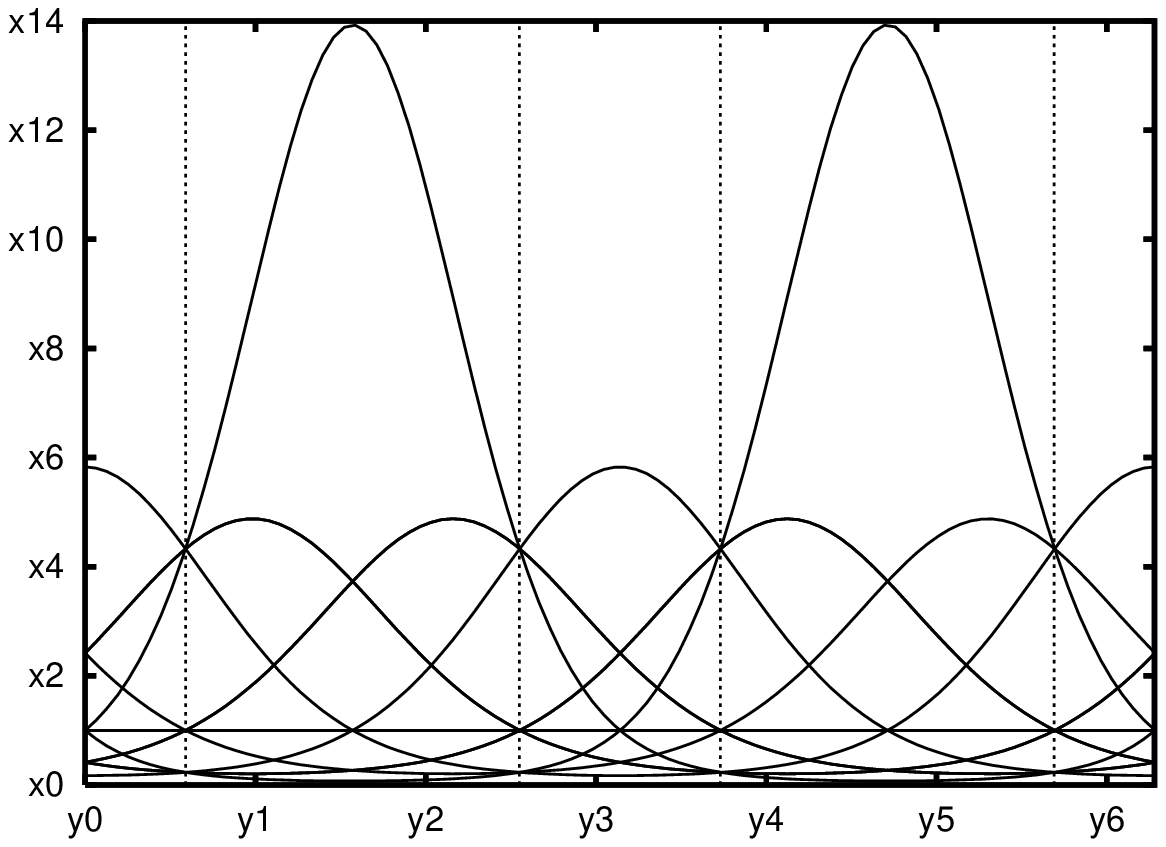}
\caption{$\overline{\Omega}_\alpha$ for Example
\ref{(1+sqrt(2))_(2+sqrt(3))_example}.}
\label{(1+sqrt(2))_(2+sqrt(3))_figure}
\end{figure}
The complete set of directional pole and zero
data~$\overline{\Omega}_\alpha$ is shown in
Figure~\ref{(1+sqrt(2))_(2+sqrt(3))_figure}. The non-expansive
directions are marked with vertical dashed lines.
\end{eg}

\begin{eg}\label{2_3_5_example}
Consider the~$\mathbb{Z}^3$-action~$\alpha$ which is dual to the
three commuting automorphisms~$\times 2,\times 3$ and~$\times 5$
on~$\mathbb{Z}[\frac{1}{30}]$. This is the dynamical system
corresponding to the~$R_3$-module
\[
M=R_3/(u_1-2, u_2-3, u_3-5).
\]
The set~$\mathcal{V}(M)\cup\mathcal{W}(M)$ consists of~$4$
characters induced by the~$2$-adic,~$3$-adic,~$5$-adic and archimedean valuations
on~$\mathbb{Q}$. Thus~$\nes_1(\alpha)$ may be obtained from the
solutions to the equation
\[
\chi^*(\mathbf{n})=1
\]
where~$\mathbf{n}\in\mathbb{R}^3$,
and~$\chi\in\mathcal{V}(M)\cup\mathcal{W}(M)$. Projecting points in
the~$4$ resulting planes to~$\sph_2$ presents~$\nes_1(\alpha)$
as a union of~$4$ great circles.

\begin{figure}[!h]
\centering
\psfrag{x0}{}
\psfrag{x1}{$1$}
\psfrag{x2}{$2$}
\psfrag{x3}{$3$}
\psfrag{x4}{$4$}
\psfrag{x5}{$5$}
\psfrag{x6}{$6$}
\psfrag{y0}{$0$}
\psfrag{y1}{$1$}
\psfrag{y2}{$2$}
\psfrag{y3}{$3$}
\psfrag{y4}{$4$}
\psfrag{y5}{$5$}
\psfrag{y6}{$6$}
\psfrag{z0}{$0$}
\psfrag{z0.5}{}
\psfrag{z1}{$1$}
\psfrag{z1.5}{}
\psfrag{z2}{$2$}
\psfrag{z2.5}{}
\psfrag{z3}{$3$}
\psfrag{z3.5}{}
\psfrag{z4}{$4$}
\includegraphics[scale=0.85]{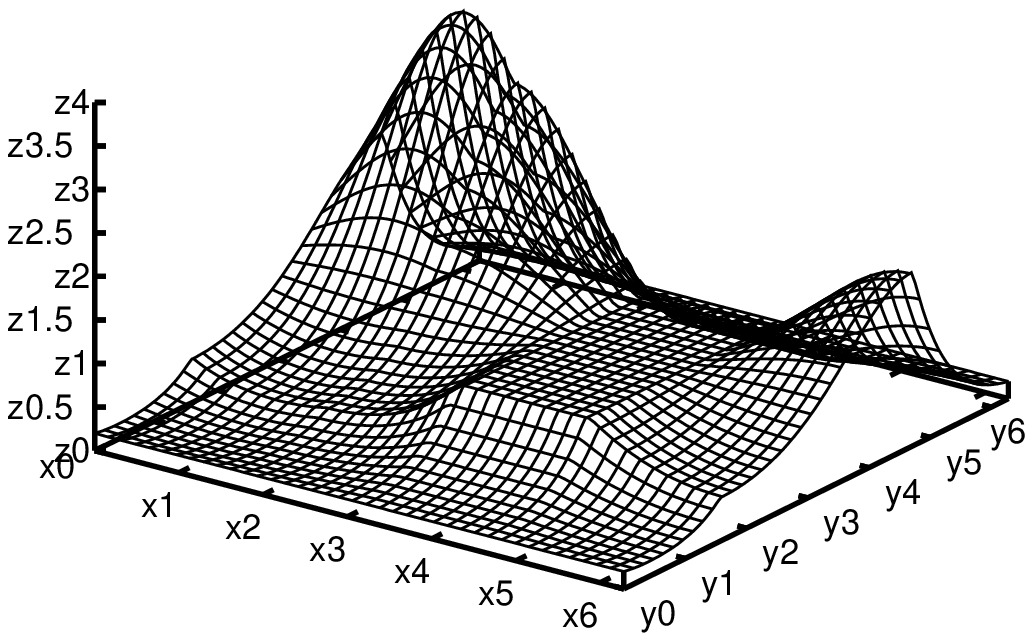}
\caption{A layer of~$\overline{\Omega}_\alpha$ for Example
\ref{2_3_5_example}.} \label{2_3_5_surface_figure}
\end{figure}

The set of directional pole and zero data~$\Omega_\alpha$ also
contains this information. Proceeding in a similar way to
Example~\ref{2_3_example},~$\overline{\Omega}_\alpha$ comprises~$2$
continuous surfaces, with curves of non-differentiability detecting
points in~$\nesn_1(\alpha)$. In terms of spherical
coordinates~$\theta,\phi\in[0,2\pi)$, these are given by
\[
\theta,\phi\in\{k\pi/2\mid k=0,1,2,3\}.
\]
Figure~\ref{2_3_5_surface_figure} shows one of the two surfaces as a
subset of~$[0,2\pi)^2\times\mathbb{R}$. The other surface intersects
this one at points detecting~$\nesv_1(\alpha)$, shown as a subset
of~$[0,2\pi)^2$ in Figure~\ref{2_3_5_intersection_figure}. This is
the locus
\[
2^{\cos\theta\sin\phi}3^{\sin\theta\sin\phi}5^{\cos\phi}=1.
\]
\begin{figure}[!h]
\centering
\psfrag{x0}{$0$}
\psfrag{x1}{$1$}
\psfrag{x2}{$2$}
\psfrag{x3}{$3$}
\psfrag{x4}{$4$}
\psfrag{x5}{$5$}
\psfrag{x6}{$6$}
\psfrag{y0}{$0$}
\psfrag{y1}{$1$}
\psfrag{y2}{$2$}
\psfrag{y3}{$3$}
\psfrag{y4}{$4$}
\psfrag{y5}{$5$}
\psfrag{y6}{$6$}
\includegraphics[scale=0.85]{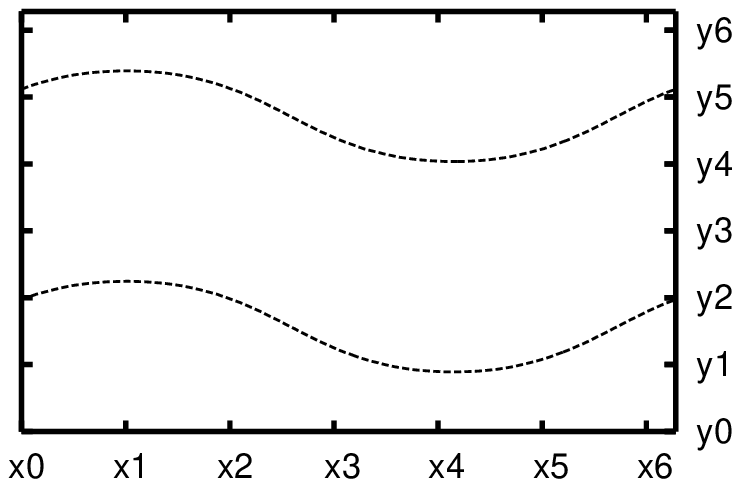}
\caption{$\nesv_1(\alpha)$ for Example \ref{2_3_5_example}.}
\label{2_3_5_intersection_figure}
\end{figure}
\end{eg}

\begin{eg}
A very important feature of algebraic dynamical systems in higher
rank is the possibility of \emph{genuinely partially hyperbolic
systems.} Damjanovi{\'c} and Katok~\cite{MR2119035} show that
for~$r\ge 6$ there are such actions on the~$r$-torus. The following
example is taken from~\cite{dkpreprint} and its directional zeta
functions are discussed in~\cite[Ex.~4.4]{miles_zeta}, where the
following description may be found.
Let
\[
f(x)=x^6-2x^5-5x^4-3x^3-5x^2-2x+1
\]
and write~$k=\mathbb Q(a)$, where~$a$ is a complex root of~$f$.
Then~$a$ and
\[
b=2a^5-6a^4-3a^3-6a^2-6a
\]
are fundamental units for the ring of integers in~$k$,
and~$M=\mathbb Z[a]$ is an~$R_2$-module under the
substitutions~$u_1=a$,~$u_2=b$. Since
\[
\{1, a, a^2, a^3, a^4, a^5\}
\]
is an integral basis for~$M$,~$\alpha_M$ is a~$\mathbb Z^2$-action
on~$X_M=\mathbb T^6$, the~$6$-torus. Moreover,~$M$ is a domain and
so it has one associated prime ideal~$\mathfrak p$ in~$R_2$ (this
ideal is the kernel of the substitution map), and~$\coht(\mathfrak
p)=1$. It follows that~$\alpha_M$ has entropy rank one. The only
places unbounded on~$M=R_2/\mathfrak p$ are archimedean; one of
these is complex and four are real. The complex place~$w$
has~$|a|_w=|b|_w=1$ (this is the part of the spectrum of the action
that is not hyperbolic) and so~$\alpha_M^{\mathbf n}$ is
non-expansive for all~$\mathbf n\in\mathbb Z^2$. The multiplicative
independence of~$a$ and~$b$ implies that~$\alpha_M^{\mathbf n}$ is
ergodic for every~$\mathbf n\in\mathbb Z^2$.
Thus~\cite[Lem~4.1]{miles_zeta} shows that~$\zeta_{\mathbf
n}$ is rational for every~$\mathbf n\in\mathbb Z^2$. This
action is not expansive, so the results above do not apply -- in
particular,~$\nes_1(\alpha)$ consists of all directions. The
crossing set and the non-smooth set for~$\overline{\Omega}_\alpha$
may be found (in principle -- there are a great many functions
involved) and these points correspond to directions~$\mathbf
n\in\mathbb R^2$ with the property that the product of some subset
of the set of eigenvalues of the integer matrix defining the
automorphism~$\widehat{\alpha_M^{\mathbf n}}$ is~$1$. Thus the
portrait obtained from~$\overline{\Omega}_\alpha$ contains those
directions in which the action is non-expansive transverse to the
two-dimensional stable foliation arising from the common
two-dimensional eigenspace for the two commuting matrices defining
the automorphisms~$\widehat{\alpha_M^{(1,0)}}$
and~$\widehat{\alpha_M^{(0,1)}}$ on which they act like a pair of
irrational rotations.
\end{eg}

\begin{eg}
Our emphasis has been on~$\mathbb Z^d$-actions, but the notion of
non-expansive subdynamics makes sense for actions of countable
abelian groups. Miles~\cite{milesmonatshefte} has extended Schmidt's
characterization~\cite{MR1069512} of expansive~$\mathbb Z^d$-actions
to this setting. When the acting group has infinite torsion-free
rank its integral group ring is no longer Noetherian as a ring, and
new phenomena are possible. The simplest example in this setting
shows that Theorem~\ref{main_theorem} does not extend without
modification.

Consider the natural action~$\alpha$ of~$\mathbb Q^{\times}_{>0}$
on~$X=\widehat{\mathbb Q}$ dual to the action~$x\overset{\times
r}{\mapsto}rx$ of~$r\in\mathbb Q^{\times}_{>0}$ on~$\mathbb Q$.
Algebraically, this is associated to a cyclic (hence Noetherian) module over the
integral group ring~$\mathbb Z[\mathbb Q^{\times}_{>0}]$, but that
is deceptive since the ring is not Noetherian. A neighborhood of the
identity in~$X$ is isometric to an open set in the adele
ring~$\mathbb Q_{\mathbb A}$ (see~\cite{MR961739} for the details);
in this picture any~$x\in X\setminus\{0\}$ must have some real
or~$p$-adic coordinate not equal to zero. Multiplication by a large
real rational, or a large negative power of that prime~$p$, will
move such a point far from the identity, showing the whole action to
be expansive. Just as in Boyle and Lind~\cite{MR1355295}, any
non-expansive subgroup lies in a non-expansive subgroup~$\Gamma$
with~$\mathbb Q^{\times}_{>0}/\Gamma\cong\mathbb Z$ (that is, lies
in some subgroup in what one is tempted to
call~$\nes_{\infty-1}(\alpha)$). Any such subgroup acts
non-expansively, since it must omit some prime~$p$, and will
therefore act isometrically on that direction in~$X$.

On the other hand, it is clear that the
set~$\overline{\Omega}_\alpha$ does not detect this rich
non-expansive structure. Since~$\widehat{X}$ is a field, every
map~$\alpha^r$ for~$r\in\mathbb Q^{\times}_{>0}$
has~$\zeta_{\alpha^r}(z)=\frac{1}{1-z}$, so the
set~$\overline{\Omega}_\alpha$ comprises the graph of the constant
function~$1$.
\end{eg}

\bibliographystyle{amsplain}
\bibliography{poles_and_zeros}

\end{document}